\documentclass[12pt,reqno]{amsart}
\usepackage{a4}
\usepackage{graphicx}

\newtheorem{theorem}{Theorem}
\newtheorem{lemma}{Lemma}

\def\minor{\mathop{\rm minor}\nolimits}
\def\Kernel{\mathop{\rm Ker}\nolimits}
\def\Image{\mathop{\rm Im}\nolimits}

\author{Igor G. Korepanov}
\title[Invariants of three-manifolds from four-dimensional geometry]{Invariants of three-dimensional manifolds from four-dimensional Euclidean geometry}
\begin{document}

\begin{abstract}
This is the first in a series of papers where we will derive invariants of three-manifolds and framed knots in them from the geometry of a manifold pseudotriangulation put in some way in a {\em four\/}-dimensional Euclidean space. Thus, the elements of the pseudotriangulation acquire Euclidean geometric values such as volumes of different dimensions and various kinds of angles. Then we construct an acyclic complex made of differentials of these geometric values, and its torsion will lead, depending on the specific kind of this complex, to some manifold or knot invariants. In this paper, we limit ourselves to constructing a simplest kind of acyclic complex, from which a three-manifold invariant can be obtained.
\end{abstract}

\maketitle

\section{Introduction}

Manifold and knot invariants can be built on the basis of some kind of geometrization of triangulation simplexes. The method for their construction, in its current state, looks as follows: an acyclic complex of vector spaces and their linear mappings is built, where the vector spaces consist of differentials of geometric quantities ascribed to elements of triangulation, such as Euclidean lengths of edges. Then the {\em torsion\/} of this complex is studied. Miraculously, it turns out that it behaves in a beautiful way under {\em Pachner moves\/} --- elementary rebuildings of a triangulation of a given manifold.

Different versions of this construction have been carried out, in full or in part, for three-dimensional (see references below in this paragraph) and four-dimensional~\cite{33,24,15,sigma} manifolds. There is very little doubt that it will work for manifolds of greater dimensions as well. Returning to three-dimensional manifolds, the most straightforward geometrization for their triangulations seems to be one using Euclidean three-dimensional geometry. Such construction (although without explicit mentioning of acyclic complexes) was proposed in paper~\cite{Korepanov-JNMP}, and further development {\em for three-dimensional manifolds\/} went along the following lines:
\begin{itemize}
\item modification of the construction making use of nontrivial representations of the manifold's fundamental group or knot group~\cite{KM,Martyushev1,Martyushev2,Martyushev3,Korepanov-FPM},
\item modification for {\em framed\/} knots in manifolds, not using the mentioned nontrivial representations~\cite{DKM},
\item modification using {\em spherical\/} (but still three-dimensional) geometry instead of Euclidean~\cite{TW},
\item construction of a {\em state-sum model\/} using Euclidean geometry~\cite{BF},
\item the use (somewhat paradoxical) of {\em two-dimensional\/} affine area-preserving geometry~\cite{KM-SL2,K-SL2} instead of Euclidean three-dimensional (the manifold, however, remains three-dimensional!).
\end{itemize}

When working on this subject, we always have in mind our main goal --- construction of {\em quantum\/} invariants in {\em greater\/} dimensions. We even hope that efficient constructions will be developed which will find applications in mathematical physics as well. One way for constructing quantum invariants, formulated in paper~\cite{sigma}, was suggested by the following way of quantizing {\em integrable models in mathematical physics\/}:
$$
\begin{matrix}\textrm{classical}\\ \textrm{scalar model}\end{matrix}
\quad \stackrel{\rm generalization}{\longrightarrow} \quad
\begin{matrix}\textrm{classical}\\ \textrm{multicomponent}\\ \textrm{model}\end{matrix}
\quad \stackrel{\rm reduction}{\longrightarrow} \quad
\begin{matrix}\textrm{quantum}\\ \textrm{model}\end{matrix}\;.
$$

In this paper we, limiting ourselves, for the time being, to three-dimensional manifolds, make a step on the way to {\em multicomponentness}. This step looks modest: instead of three-dimensional, we are using {\em four-dimensional\/} Euclidean geometry. The idea is, however, clear: to pass on, in future, from three-dimensional manifolds to manifolds of any dimensions, and from four-dimensional geometry --- to some very high-dimensional geometry, from which, we hope, we will be able to fish out something {\em quantum\/} by means of a proper reduction.

\section{The main algebraic complex}

\subsection{Pseudotriangulated manifold and its geometrization}
\label{ss_geom}

Let $M$ be a pseudotriangulated, compact, connected, orientable three-dimen\-sional manifold. {\em Pseudo\/}tri\-an\-gu\-la\-tions differ from tringulations in the proper sense, or combinatorial triangulations, in that in a pseudotriangulation, the boundary of a simplex may contain a simplex of a smaller dimension several times, and a simplex is {\em not}, in general, determined uniquely by the set of its vertices. Pseudotriangulations contain in many cases a smaller number of elements (vertices, egdes, two-dimensional faces and tetrahedra) than combinatorial triangulations, and so are more convenient for calculations. Nevertheless, we impose one restriction on our pseudotriangulations: we require that, for any tetrahedron, {\em all four\/} its vertices must be {\em different}. When we consider sequences of simplicial moves (local rebuildings) of a pseudotriangulation, this requirement must be obeyed at every step. 

In this paper, we also require that $M$ must be {\em without a boundary}, although we plan to apply our constructions to manifolds with boundaries in forthcoming papers.

We also fix once and for all an orientation of~$M$. In terms of a pseudotriangulation, this means that every tetrahedron becomes oriented, that is, an ordering of its vertices up to an {\em even\/} permutation is given, and the orientations of different tetrahedra are {\em consistent\/}: for instance, two adjacent tetrahedra $EABC$ and $ABCD$ have consistent orientations. If the contrary is not stated explicitly, below we consider all tetrahedra only with their ``right'' orientation.

We geometrize our pseudotriangulation in the following way: we ascribe to its every vertex four Euclidean coordinates. The coordinates of a vertex~$A$ will be denoted $(w_A,x_A,y_A,z_A)\in \mathbb R^4$. These coordinates are arbitrary with the only restriction: no five different vertices can fall into the same affine hyperlane in~$\mathbb R^4$. In particular, if our pseudotriangulation is a triangulation in the proper sense, we get thus a {\em piecewise-linear imbedding\/} of~$M$ into~$\mathbb R^4$, i.e., a locally one-to-one piecewise-linear mapping. The space~$\mathbb R^4$ will be understood as having a fixed orientation given by the order of its axes $w,x,y$, and~$z$; we call this orientation {\em positive}.

Given the coordinates of vertices, we get in a natural way other geometric values for the elements of pseudotriangulation, such as volumes of simplices of different dimensions and different sorts of Euclidean angles. We introduce some notations for them. The length of an edge~$AB$ will be denoted~$l_{AB}$. The area of a two-dimensional face $ABC$ will be denoted~$S_{ABC}$. The volume of a tetrahedron $ABCD$ will be denoted~$V_{ABCD}$. All these values are supposed to be {\em positive}. In contrast with this, when we need a {\em four\/}-volume $\mathcal V_{ABCDE}$ of a 4-simplex with vertices $A,B,C,D,E$, this will be the {\em oriented\/} volume:
$$
\mathcal V_{ABCDE} \stackrel{\rm def}{=} \frac{1}{24}\det\left(
\begin{array}{c}
\overrightarrow{AB}\\ \overrightarrow{AC}\\ \overrightarrow{AD}\\ \overrightarrow{AE}
\end{array}
\right),
$$
where $\overrightarrow{AB}=(w_B-w_A,x_B-x_A,y_B-y_A,z_B-z_A)$, etc.

Next we introduce two sorts of angles. First, these are {\em inner\/} dihedral angles in a tetrahedron; we denote them $\varphi$, maybe with subscripts whose meaning is specified when necessary. The angles~$\varphi$ are supposed to take values in the interval $(0,\pi)$. Second, we will need angles between {\em adjoint tetrahedra}, defined as follows. Let there be two tetrahedra $EABC$ and~$ABCD$ having the common face~$ABC$. By definition, the angle~$\vartheta_{ABC}$ between them is the {\em exterior\/} dihedral angle at the two-face $ABC$ in the four-simplex $ABCDE$, which is, moreover, taken with the minus sign in case if $\mathcal V_{ABCDE}<0$. Note that
$$
\vartheta_{ABC}=\vartheta_{BAC},
$$
because this change of the orientation of face~$ABC$ implies also the change of order in which the oriented tetrahedra $EABC$ and $ABCD$ come.

\subsection{Three-component deficit angles}
\label{ss_da}

Three-component deficit angles are a key notion in this paper. Let there be an edge~$BC$ in the pseudotriangulation of~$M$, and let its link consist of edges $A_1A_2,\ldots,A_nA_1$ in such way that tetrahedra $A_1A_2BC,\ldots,\allowbreak A_nA_1BC$ have the right orientation. These tetrahedra form together the {\em star\/} of~$BC$.

Now we geometrize this star in a way slightly different from Subsection~\ref{ss_geom}. We begin with ascribing lengths to all edges in the star in such way that $A_1A_2BC,\allowbreak \ldots,\allowbreak A_nA_1BC$ become Euclidean tetrahedra. Then we ascribe some real numbers $\vartheta_{A_1BC},\ldots,\allowbreak \vartheta_{A_nBC}$, defined modulo $2\pi$, to the respective faces $A_1BC,\ldots,\allowbreak A_nBC$. {\em After\/} this, we start trying to put this structure in~$\mathbb R^4$ in the following way. First, we put in~$\mathbb R^4$ the tetrahedron $A_1A_2BC$, i.e., give Euclidean coordinates $w,x,y,z$ to its vertices in such way that distances between vertices coincide with the given lengths of edges. We can say that we {\em associate a coordinate system\/} with tetrahedron~$A_1A_2BC$. Then, given the number~$\vartheta_{A_2BC}$, we can give Euclidean coordinates to vertex~$A_3$, in a unique way provided we want to obtain again the right edge lengths and $\vartheta_{A_2BC}$ as the angle between tetrahedra $A_1A_2BC$ and $A_2A_3BC$ in the sense of Subsection~\ref{ss_geom}. We can say that we {\em extend the coordinate system\/} to tetrahedron~$A_2A_3BC$. Continuing this way and making one full revolution around~$BC$, we get as a result a {\em new\/} coordinate system for~$A_1A_2BC$, with new axes $w^{\rm new},x^{\rm new},y^{\rm new},z^{\rm new}$.

The transition from old to new coordinate system defines an element of the group of Euclidean motions of~$\mathbb R^4$. However, as old and new coordinates obviously coincide for vertices $B$ and~$C$, we actually get a rotation within the space orthogonal to~$BC$, i.e., an element of group~${\rm SO}(3)$, which we denote as~$\omega_{BC}$ and call {\em three-component deficit angle}, or {\em dicrete curvature}, around~$BC$. The case $\omega_{BC}=1\in {\rm SO}(3)$ corresponds exactly to the situation where the whole star of~$BC$ can be put in~$\mathbb R^4$ without ``cracks''.

We do not specify here how we represent the values $\omega$ in matrix form. This is because we will mostly need their {\em infinitesimal\/} versions $d\omega$, for $\omega$ infinitely close to unity, and we will explain in Subsection~\ref{ss_seq} how we represent infinitesimal quantities.

\subsection{A sequence of linear spaces and linear mappings}
\label{ss_seq}

We will now construct our ``main'' algebraic complex, consisting of based linear spaces of differentials of geometric values and linear mappings between them. The word ``based'' means that a vector space is taken together with some basis in it chosen in some natural way. This means that elements of the vector space can be represented as {\em column vectors}, and linear mappings as {\em matrices}.

We call our complex ``main'' because, besides studying it as it is, we will also be interested in its various modifications. It is expected that this complex and/or its modifications are {\em acyclic\/} in many interesting cases; although we do not prove general theorems about acyclicity in this paper, we will show it at least on examples in the next paper(s) in this series. For an acyclic complex, one can calculate its {\em torsion}, and we show in Section~\ref{s_moves} how to obtain a topological invariant from this.

Consider the following sequence of vector spaces and linear mappings. We first just write it out, and then we explain in this Subsection the used notations and give the definitions for both spaces and mappings. The proof that our sequence is actually an algebraic complex is given in Subsection~\ref{ss_c}. So, our sequence is:
\begin{equation}
0\rightarrow \mathfrak e_4
\stackrel{f_1}{\rightarrow} (dx)
\stackrel{f_2}{\rightarrow} \left( \begin{array}{c} dl\\ \oplus\\ d\vartheta \end{array} \right)
\stackrel{f_3}{\rightarrow} (d\omega)
\stackrel{f_4}{\rightarrow} (d\rho)
\stackrel{f_5}{\rightarrow} \left( \begin{array}{c} d\alpha\\ \oplus\\ d\beta \end{array} \right)
\rightarrow 0\,.
\label{mc}
\end{equation}

We assume that vertices of the pseudotriangulation of~$M$ have been put in~$\mathbb R^4$ as described in Subsection~\ref{ss_geom}. The vector space~$\mathfrak e_4$ is the Lie algebra of infinitesimal Euclidean motions of~$\mathbb R^4$. The vector space denoted~$(dx)$ consists of column vectors of differentials $dw_A,dx_A,dy_A,dz_A$ of all coordinates for all vertices~$A$ in the pseudotriangulation. The definition of mapping~$f_1$ is obvious: it reflects how an infinitesimal motion of~$\mathbb R^4$ changes the coordinates $(w_A,x_A,y_A,z_A)$ of all points~$A$. Explicitly, we will write it as follows. We identify algebra~$\mathfrak e_4$ with the vector space of columns
\begin{equation}
(dr_1,dr_2,dr_3,dr_4,a_{12},a_{13},a_{14},a_{23},a_{24},a_{34})^{\rm T},
\label{e4}
\end{equation}
the superscript $\rm T$ standing for matrix transposing. Vector~(\ref{e4}) gives, by definition, the following differentials of coordinates of a given vertex~$A$ due to mapping~$f_1$:
\begin{equation}
\left( \begin{array}{c} dw_A\\dx_A\\dy_A\\dz_A \end{array} \right) = \left( \begin{array}{c} dr_1\\dr_2\\dr_3\\dr_4 \end{array} \right)
+ \mathcal A \left( \begin{array}{c} w_A\\x_A\\y_A\\z_A \end{array} \right),
\label{Acal}
\end{equation}
where $\mathcal A$ is the {\em antisymmetric\/} matrix with elements~$a_{ij}$. So, the vector~$d\vec r=(dr_1,dr_2,\allowbreak dr_3,dr_4)^{\rm T}$ represents a translation, while $\mathcal A$ is an element of Lie algebra~$\mathfrak s \mathfrak o (4)$.

The third vector space, $(dl\oplus d\vartheta)$, is a direct sum of the space~$(dl)$ consisting of differentials of all edge lengths in the pseudotriangulation, and the space~$(d\vartheta)$ consisting of differentials of values~$\vartheta$ (see Subsection~\ref{ss_geom}) for all two-dimensional faces. The definition of mapping~$f_2$ is again obvious: if the coordinates of all vertices are given, they determine naturally all $l$ and~$\vartheta$. Mapping~$f_2$ is the infinitesimal version of this correspondence.

We have thus described $f_2$ geometrically, without giving explicit formulas like~(\ref{Acal}). We will also give in this Section only a geometric description of mapping~$f_3$. Nevertheless, we will provide explicit formulas when we need them in our calculations in Section~\ref{s_moves}.

The fourth vector space, denoted $(d\omega)$, consists of infinitesimal three-compo\-nent deficit angles for all edges in the pseudotriangulation. An infinitesimal three-compo\-nent deficit angle~$d\omega_{BC}$ is by definition an infinitesimal rotation around the edge~$BC$, thus represented by an antisymmetric $4\times 4$ matrix~$\mathcal B$ with the property $\mathcal B\,\overrightarrow {BC}=0$. We require that
\begin{equation}
d\omega_{CB}=-d\omega_{BC},
\label{cb}
\end{equation}
this antisymmetry reflects the fact that the direction in which we go around $BC$ is determined by the direction of~$BC$ according to the corkscrew rule, and is in agreement with the definition of mapping~$f_3$ which we will give soon.

On the other hand, we want to think of the space $(d\omega)$ as consisting of column vectors of dimension~$3N_1$, where $N_1$ is the number of edges in the pseudotriangulation. To each edge belong three components, which we define in the following way. Choose some positively oriented Euclidean coordinate system with axes $w',\allowbreak x',\allowbreak y',\allowbreak z'$ in such way that axis~$w'$ goes in the direction of vector~$\overrightarrow {BC}$. In this coordinate system, $\mathcal B$ is transformed to a matrix with elements~$(b_{ij}')$ such that $b_{12}'=b_{13}'=b_{14}'=0$. And the three elements $b_{34}',b_{42}',b_{23}'$ are by definition the three mentioned components of the column vector.

The mapping~$f_3$ shows which deficit angles result from given deformations of lengths~$l$ and angles~$\vartheta$ and goes as follows. Imagine that initially all lengths~$l$ and angles~$\vartheta$ come from some vertex coordinates due to usual Euclidean formulae. If we then change $l$'s and $\vartheta$'s slightly but otherwise arbitrarily, the whole structure can no longer be placed in~$\mathbb R^4$. Still, we can define deficit angles as in Subsection~\ref{ss_da}. If now the changes of $l$'s and $\vartheta$'s are {\em infinitesimal}, we get infinitesimal three-compo\-nent deficit angles~$d\omega$ which are nothing but infinitesimal rotations around corresponding edges. It is not hard to see that matrices representing these rotations {\em do not depend\/} on the tetrahedron from which we start to drag a coordinate system around an edge (in Subsection~\ref{ss_da}, the edge was called~$BC$, and the tetrahedron $A_1A_2BC$). So, mapping~$f_3$ by definition provides exactly this way of obtaining $d\omega$'s from $dl$'s and~$d\vartheta$'s.

The fifth vector space, denoted $(d\rho)$, consists of antisymmetric $4\times 4$ matrices $d\rho_A$ corresponding to every vertex~$A$. To represent $d\rho_A$ as a column vector, we act again in the same way as for matrix~$\mathcal A$ in~(\ref{Acal}). The mapping~$f_4$ gives, by definition, the following~$d\rho_A$ from given~$d\omega$'s:
\begin{equation}
d\rho_A = \sum_B d\omega_{AB},
\label{rho}
\end{equation}
where the sum is taken over all vertices~$B$ joined to~$A$ by edges. In (\ref{rho}), both $d\rho_A$ and $d\omega_{AB}$ are considered, of course, as $4\times 4$ matrices --- elements of the same algebra~$\mathfrak s \mathfrak o (4)$.

The last, sixth vector space, denoted $(d\alpha \oplus d\beta)$, is ten-dimen\-sional, like the first space~$\mathfrak e_4$. Here $d\alpha$ is {\em one\/} matrix --- element of~$\mathfrak s \mathfrak o (4)$, and $d\beta$ is one four-dimen\-sional vector. By definition, mapping~$f_5$ makes the following $d\alpha$ and~$d\beta$ from given~$d\rho_A$:
\begin{equation}
d\alpha = \sum_{{\rm over\;all\;}A} d\rho_A,\qquad d\beta = \sum_{{\rm over\;all\;}A} d\rho_A\, \vec r_A,
\label{alpha-beta}
\end{equation}
where $\vec r_A=(w_A,x_A,y_A,z_A)^{\rm T}$ is the radius vector of vertex~$A$, and $d\rho_A\, \vec r_A$ is the result of the action of the linear operator~$d\rho_A$ on this vector. The sums in~(\ref{alpha-beta}) are taken over all vertices~$A$ in the pseudotriangulation.

\subsection{Sequence (\ref{mc}) is an algebraic complex}
\label{ss_c}

\begin{theorem}
\label{th_c}
Sequence (\ref{mc}) is an algebraic complex, i.e., all the compositions of neighboring linear operators $f_2f_1$, $f_3f_2$, $f_4f_3$ and $f_5f_4$ are zero operators.
\end{theorem}

\begin{proof}
The proof of Theorem~\ref{th_c} follows from simple geometric considerations.

(a) $f_2f_1=0$ is just a differential version of the fact that Euclidean motions do not change lengths and angles.

(b) $f_3f_2=0$ is a differential version of the fact that, if lengths and angles for our pseudotriangulation are found from given coordinates of vertices, then there are no deficit angles (i.e., all $\omega\equiv 1$ and $d\omega=0$).

(c) $f_4f_3=0$ can be obtained as follows. Let the edges going out of a vertex~$A$ be $AB_1,\ldots,\allowbreak AB_n$. Suppose our geometrization is now made in the sense of Subsection~\ref{ss_da}: edge lengths~$l$ and angles~$\vartheta$ between adjacent tetrahedra are given, but no vertex coordinates. When we go around edge~$AB_i$ in manifold~$M$ in the same way as we did it for edge~$BC$ in Subsection~\ref{ss_da}, we get a rotation of a Euclidean coordinate system associated with the tetrahedron from which we have started. Likewise, we get such a rotation if we go along some closed trajectory starting within a fixed tetrahedron $AB_1B_2B_3$ and going always through tetrahedra belonging to the star of vertex~$A$, with the condition that this trajectory is allowed to intersect two-dimensional faces but not edges. It is obvious that if one such trajectory can be continuously deformed into another in such way that it never intersects any edge, the resulting rotation --- element of group~${\rm SO}(4)$ --- will be the same for both trajectories.

Consider now the infinitesimal version of this situation, where all rotations are infinitely close to unity, differing from it by some elements $d\omega\in\mathfrak s\mathfrak o(4)$. First, it is not hard to see that if we go from a given tetrahedron $AB_1B_2B_3$ along a trajectory surrounding just {\em one\/} edge~$AB_i$, where $B_i$ does not need to coincide with any of $B_1,B_2,B_3$, then the corresponding $d\omega=d\omega_i$ does {\em not\/} depend on the exact way used to reach~$AB_i$ and go back. Second, if a trajectory surrounds several edges, then $d\omega$ is just a sum of corresponding~$d\omega_i$. Third, given that the link of~$A$ is a sphere, the trajectory surrounding {\em all\/} edges is the same as the trajectory surrounding {\em no\/} edges. This means exactly that the r.h.s.\ of~(\ref{rho}) vanishes if the $d\omega$'s come from some $dl$'s and~$d\vartheta$'s.

(d) $f_5f_4=0$. First, each $d\omega$, for instance, $d\omega_{AB}$, yields due to~$f_4$ mutually opposite values~$d\rho$ on two ends of the edge: 
$$
d\rho_A=-d\rho_B\quad \hbox{ if only } d\omega_{AB}\ne 0 \hbox{ of all } d\omega,
$$
see (\ref{cb}) and (\ref{rho}). Now it is clear that the first equation~(\ref{alpha-beta}) gives zero if the $d\rho$ come from {\em any\/}~$d\omega$.

Second, considering in the same way only one nonzero $d\omega=d\omega_{AB}$, we get from the second equation~(\ref{alpha-beta})
$$
d\beta=d\rho_A\,\vec r_A+d\rho_B\,\vec r_B=d\rho_A(\vec r_A-\vec r_B)=d\omega_{AB} \,\overrightarrow{BA}=0,
$$
because a rotation around edge~$AB$ leaves this edge intact.
\end{proof}

\section{Moves $2\leftrightarrow 3$ and $1\leftrightarrow 4$ and their invariant}
\label{s_moves}

\subsection{Piecewise-linear manifolds and Pachner moves}

A very good introduction in triangulated manifolds and Pachner moves is paper~\cite{L}. Although it deals with combinatorial triangulations, it is just a small exercise to extend its results onto our situation with pseudotriangulated manifolds subject to conditions stated in the beginning of Subsection~\ref{ss_geom}. What matters for us is that any (pseudo)tri\-an\-gu\-la\-tion of our manifold~$M$ can be transformed into any other (pseudo)tri\-an\-gu\-la\-tion by a finite sequence of Pachner moves, which are in three dimensions moves $2\to 3$ and $1\to 4$, explained in Subsections \ref{ss_2-3} and~\ref{ss_1-4} respectively, and the inverses of these moves. Thus, a quantity invariant under all Pachner moves is an invariant of a piecewise-linear manifold. Given the well-known fact that, in three dimensions, the categories of piecewise-linear, topological and differentiable manifolds essentially coincide, our invariant will be an invariant of a topological or differentiable manifold as well.

\subsection{Acyclicity and torsion of complex (\ref{mc})}

We will build an invariant of Pachner moves on the base of the torsion of complex~(\ref{mc}) {\em assuming its acyclicity}. Acyclicity, or exactness, for a chain complex like~(\ref{mc}) means that the image of any homomorphism (in our case --- linear mapping of vector spaces) coincides exactly with the kernel of the next homomorphism:
$$
\Kernel f_1=0,\quad \Image f_1=\Kernel f_2,\quad \ldots ,\quad \Image f_5=(d\alpha\oplus d\beta).
$$
For an acyclic complex~(\ref{mc}), we define its {\em torsion\/} as
\begin{equation}
\tau=\frac{\minor f_2\cdot \minor f_4}{\minor f_1\cdot \minor f_3\cdot \minor f_5},
\label{tau}
\end{equation}
where we take for each matrix $f_1,\ldots,f_5$ a nonzero minor of the maximal order. Further conditions are: the basis of each of six vector spaces in~(\ref{mc}) must be completely divided between the two minors of matrices standing to the left and to the right of this vector space. For example, if the rows of $\minor f_3$ correspond to some chosen basis vectors in space~$(d\omega)$, then the columns of $\minor f_4$ must correspond to exactly all the rest of basis vectors in~$(d\omega)$.

In terms of the textbook~\cite{Turaev}, we are using the definition of torsion based on the notion of a {\em non-degenerate $\tau$-chain\/} which consists, in our case, of submatrices of matrices $f_1,\ldots,f_5$ corresponding to the mentioned minors. As is known, such torsion~$\tau$ does not depend on the specific choice of minors, except for its sign which changes under an odd permutation of basis vectors in any space. Below, we always consider torsion~$\tau$, determinants related to it and equalities contatining them {\em up to a sign}.

It is also of use to have in mind the (equivalent) definition of torsion given in Section~1 of book~\cite{Turaev}. We do not give it here, but the essential point for us is that it is obvious from this definition how to justify the orthogonal rotations of coordinate axes which we will freely use below in our proofs of Theorems \ref{th_2-3} and~\ref{th_1-4}. The interested reader will see this from the formula given in Remark~1 in Subsection~1.4 of~\cite{Turaev} and the fact that the determinant of an orthogonal rotation is unity.

So, below we always assume that complex~(\ref{mc}) is acyclic. This convention needs to be made more precise in the following way. To construct a manifold invariant, we are going to investigate how $\tau$ changes under Pachner moves $2\leftrightarrow 3$ and $1\leftrightarrow 4$. It turns out that $\tau$ is multiplied by some `local' (belonging to the clusters of tetrahedra which undergo a transformation) factors under moves $2\to 3$ and $1\to 4$; then under the reverse moves $\tau$ is, of course, divided by the same factors. It is clear from the very fact that the mentioned `local' factors are correctly defined and are neither zero nor infinity that, if the complex is acyclic before a move $2\to 3$ or $1\to 4$, then it is so after the move, and vice versa. To make this statement completely accurate, one can use Lemma~2.5 from book~\cite{Turaev}, stating that ``A chain complex over a field is acyclic if and only if it has a non-degenerate $\tau$-chain''. It will be clear from the proofs of our Theorems \ref{th_2-3} and~\ref{th_1-4} how to choose submatrices for a non-degenerate $\tau$-chain for a complex resulting from a Pachner move.

\subsection{Move \boldmath $2\to 3$}
\label{ss_2-3}

\begin{figure}
\begin{center}
\includegraphics[scale=1.5]{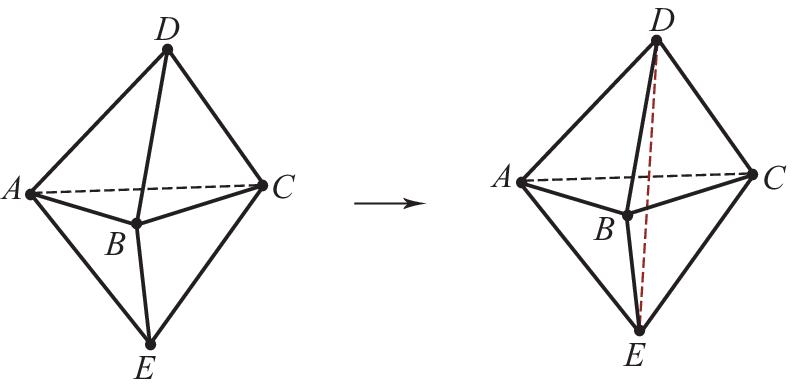}
\end{center}
\caption{Pachner move $2\to 3$}
\label{fig1}
\end{figure}

Under a move $2\to 3$, two adjacent tetrahedra $ABCD$ and $EABC$ are replaced with three tetrahedra $ABED$, $BCED$ and~$CAED$, see Figure~\ref{fig1}.
Thus, the changes do not affect the {\em vertices\/} (0-cells) of the simlicial complex and vector spaces $(dx)$ and~$(d\rho)$ corresponding to them. There is one additional edge (1-cell) $DE$ in the right-hand-side diagram of Figure~\ref{fig1}, compared to the left-hand-side diagram. Consequently, the space~$(dl)$ will acquire one additional basis vector~$dl_{DE}$, and the space~$(d\omega)$ will acquire three basis vectors --- components of~$d\omega_{DE}$. To emphasize the three-component character of this differential, we write below $d\vec\omega_{DE}$ instead of~$d\omega_{DE}$ and, generally, $d\vec\omega_i$ instead of~$d\omega_i$, where $i$ is an edge. Changes in two-dimensional faces imply that basis vector~$d\vartheta_{ABC}$ is excluded from the space~$(d\vartheta)$, while three new vectors $d\vartheta_{ADE}$, $d\vartheta_{BDE}$ and~$d\vartheta_{CDE}$ are added to it.

\begin{theorem}
\label{th_2-3}
Under the move $2\to 3$ depicted in Figure~\ref{fig1}, torsion~$\tau$ of complex~(\ref{mc}) is multiplied by
\begin{equation}
\frac{2S_{ADE}\cdot 2S_{BDE}\cdot 2S_{CDE}}{l_{DE}^3\cdot 2S_{ABC}}.
\label{tau*}
\end{equation}
\end{theorem}

\begin{proof}

First we will need the following simple lemma.
\begin{lemma}
\label{l*}
The `initial' minors in formula (\ref{tau}), i.e., those corresponding to the situation {\em before\/} move $2\to 3$, can always be chosen in such way that $\minor f_3$ contains the column corresponding to~$d\vartheta_{ABC}$.
\end{lemma}

\begin{proof}[Proof of Lemma \ref{l*}]
Consider the tetrahedra in the star of edge~$AB$. If they all lie in~$\mathbb R^4$, then $\vartheta_{ABC}$ is a function of the rest of $\vartheta_{AB\ldots}$ and the lengths of edges in these tetrahedra. Hence, the row $\frac{\partial\vartheta_{ABC}}{\partial x}$ of matrix~$f_2$ (i.e., the row corresponding to~$d\vartheta_{ABC}$) is a linear combination of its other rows. This means that a nonvanishing minor of~$f_2$ containing this row can always be replaced with a nonvanishing minor of the same size {\em not\/} containing it. Consequently, basis vector $d\vartheta_{ABC}$ can be ascribed to $\minor f_3$. This will of course cause no trouble to $\minor f_2$ because, due to the acyclicity in term~$(dl\oplus d\vartheta)$, this minor will just change in such way that ratio $\frac{\minor f_2}{\minor f_3}$ remains the same.
\end{proof}

Now we continue the proof of Theorem~\ref{th_2-3}. Let the initial $\minor f_3$ be as in Lemma~\ref{l*}. {\em After\/} the move $2\to 3$, we take a new $\minor f_3$ which contains the column~$dl_{DE}$ instead of the column~$d\vartheta_{ABC}$ in the old minor, as well as three new columns $d\vartheta_{ADE},d\vartheta_{BDE},d\vartheta_{CDE}$ and three new rows for components of~$d\vec\omega_{DE}$. So, the minors of $f_1,f_2,f_4$ and~$f_5$ remain the same, and all the change of $\tau$ is concentrated in $\minor f_3$.

We examine how the `new' matrix $f_3$, i.e., matrix corresponding to the right-hand-side diagram of Figure~\ref{fig1}, is obtained from the `old' matrix. As we remember, matrix~$f_3$ consists of partial derivatives of components of quantities~$\vec\omega_i$ with respect to lengths~$l_i$ and angles~$\vartheta_a$ ($i$ and $a$ being an edge and a two-face, respectively). At this moment, we will have to consider $\vec\omega_i$, for the `new' simplicial complex, as functions of a slightly different set of independent variables. To define them, we first introduce quantities $\vartheta_{ADE}^{(0)},\vartheta_{BDE}^{(0)},\vartheta_{CDE}^{(0)}$ as the values of corresponding angles~$\vartheta$ determined from the `zero curvature' condition
\begin{equation}
\vec\omega_{DE}=1
\label{DE_flat}
\end{equation}
with given {\em other\/} $\vartheta$ and lengths~$l$ in the simplitial complex. It is enough for us that $\vartheta_{ADE}^{(0)},\vartheta_{BDE}^{(0)}$ and $\vartheta_{CDE}^{(0)}$ exist when the other $\vartheta$ and~$l$ stay close to their initial values. Condition (\ref{DE_flat}) means that the three tetrahedra in the right-hand side of Figure~\ref{fig1} can be put together in~$\mathbb R^4$. Thus, angle~$\vartheta_{ABC}$ is correctly defined is this situation (although no face~$ABC$ is present in the right-hand side diagram of Figure~\ref{fig1}).

We consider quantities $\vec\omega_i$ for the `new' simplicial complex as functions of the following independent variables: edge lengths~$l_i$ and angles~$\vartheta_a$, except for length~$l_{DE}$, instead of which we now take {\em angle\/}~$\vartheta_{ABC}$, and angles $\vartheta_{ADE},\vartheta_{BDE}$ and~$\vartheta_{CDE}$, instead of which we take
\begin{equation}
\psi_{ADE}=\vartheta_{ADE}-\vartheta_{ADE}^{(0)},\quad \psi_{BDE}=\vartheta_{BDE}-\vartheta_{BDE}^{(0)},\quad\hbox{and}\quad\psi_{CDE}=\vartheta_{CDE}-\vartheta_{CDE}^{(0)}.
\label{e*}
\end{equation}
To make this more precise, we describe the way of obtaining `usual' independent variables $l$ and~$\vartheta$ from our `new' variables ``$\vartheta_{ABC}$, three~$\psi$'s and the rest of $\vartheta$ and~$l$'': first, we ignore the three~$\psi$'s and calculate $l_{DE}$ and $\vartheta_{ADE}^{(0)},\vartheta_{BDE}^{(0)},\vartheta_{CDE}^{(0)}$ using four-dimensional Euclidean geometry, then we add the three~$\psi$'s to the corresponding angles~$\vartheta^{(0)}$ for obtaining $\vartheta_{ADE}$, $\vartheta_{BDE}$ and~$\vartheta_{CDE}$.

Denote the Jacobian matrix of partial derivatives of components of $\vec\omega_i$ with respect to our `new' variables as~$\tilde f_3$. The fact that (\ref{DE_flat}) holds when
\begin{equation}
\psi_{ADE}=\psi_{BDE}=\psi_{CDE}=0
\label{3psi}
\end{equation}
identically in other variables --- we call them `non-psi' --- means that the derivatives of $\vec\omega_i$ with respect to these non-psi variables vanish. Besides, the rest of $\vec\omega_i$, $i\ne DE$, depend on the non-psi variables, provided (\ref{3psi}) holds, in the same way as for the initial complex (before the $2\to 3$ move), because the condition (\ref{3psi}) means that both sides of Figure~\ref{fig1} can be superposed in~$\mathbb R^4$, and the immediate consequence of this is that $\vec\omega_i$, $i\ne DE$, are the same for the initial complex and one resulting from $2\to 3$ move.

Thus, $\tilde f_3$ has the following block-triangular form:
\begin{equation}
\tilde f_3=\begin{pmatrix} \partial \vec\omega_{DE} / \partial \psi && \mathbf 0 \\[1ex] \hbox{\boldmath $\ast$} && f_3^{\rm old} \end{pmatrix},
\label{e**}
\end{equation}
where $\partial \vec\omega_{DE} / \partial \psi$ is, of course,  the $3\times 3$ matrix of derivatives of $\vec\omega_{DE}$  components with respect to variables~(\ref{e*}).

Matrix $f_3$ is obtained from $\tilde f_3$ by the right multiplication by a matrix~$g$ of partial derivatives of the set of variables $(\psi_{ADE},\psi_{BDE},\psi_{CDE}; \vartheta_{ABC}; \hbox{other } \vartheta \hbox{ and } l, \hbox{ except } l_{DE})$ with respect to the set $(\vartheta_{ADE},\vartheta_{BDE},\vartheta_{CDE}; l_{DE}; \hbox{other } \vartheta \hbox{ and } l)$:
\begin{equation}
f_3=\tilde f_3 g.
\label{e***}
\end{equation}
The semicolons between groups of variables in parentheses correspond to the block decomposition of matrix~$g$ in the following formula which also reflects a block triangular form of~$g$:
\begin{equation}
g=\left( \begin{array}{c|c} \vphantom{\big|} \mathbf 1 & \hbox{\boldmath $\ast$} \\\hline \mathbf 0 &
\begin{array}{c|c} \vphantom{\Big|}\frac{\partial\vartheta_{ABC}}{\partial l_{DE}} & \hbox{\boldmath $\ast$} \\\hline \vphantom{\Big|} \mathbf 0 & \mathbf 1
\end{array}
\end{array} \right).
\label{e4*}
\end{equation}

One sees from formulas (\ref{e**}), (\ref{e***}) and~(\ref{e4*}) that the resulting change of $\minor f_3$ under move $2\to 3$ is described as follows:
$$
(\minor f_3)^{\rm new} = (\minor f_3)^{\rm old}\cdot \det\left(\frac{\partial \vec\omega_{DE}}{\partial\psi}\right)\cdot \frac{\partial \vartheta_{ABC}}{\partial l_{DE}}.
$$
The rightmost multiplier here is calculated easily:
\begin{equation}
\frac{\partial \vartheta_{ABC}}{\partial l_{DE}}=-\frac{2S_{ABC}\,l_{DE}}{24\mathcal V_{ABCDE}},
\label{e5*}
\end{equation}
see formula (15) in \cite{multidim}; the minus sign in~(\ref{e5*}) is due to $\vartheta_{ABC}$ being an {\em exterior\/} angle. 

So, what remains is to calculate $\det\left(\frac{\partial \vec\omega_{DE}}{\partial\psi}\right)$. This determinant consists of three  columns, each representing a partial derivative of three-component~$\vec\omega_{DE}$ in some~$\psi$. Consider, for instance, the derivative $\partial\vec\omega_{DE} / \partial\psi_{ADE}$. It can be represented as a vector of {\em unit length\/} in a three-dimensional space, if we represent $d\vec\omega_{DE}$ as a three-vector according to Subsection~\ref{ss_seq}, the paragraph after formula~(\ref{cb}). As for the direction of $\partial\vec\omega_{DE} / \partial\psi_{ADE}$, it is determined (up to a sign) by the fact that the infinitesimal rotation corresponding to a differential~$d\psi_{ADE}$ takes place in the two-dimensional plane orthogonal to two-face~$ADE$.

Similar statements hold for two other derivatives $\partial\vec\omega_{DE} / \partial\psi_{BDE}$ and $\partial\vec\omega_{DE} / \partial\psi_{CDE}$. To calculate the determinant made of these three unit vectors, one must know angles between them, and these angles are the same as between vectors $\overrightarrow{AD}_{\bot\overrightarrow{DE}}$, $\overrightarrow{BD}_{\bot\overrightarrow{DE}}$ and $\overrightarrow{CD}_{\bot\overrightarrow{DE}}$, where `$\bot\overrightarrow{DE}$' means the component of the respective vector {\em orthogonal to\/}~$\overrightarrow{DE}$. From these three latter vectors, one can readily make unit vectors, replacing them with
\begin{equation}
\frac{l_{DE}}{2S_{ADE}}\overrightarrow{AD}_{\bot\overrightarrow{DE}}\;,\quad \frac{l_{DE}}{2S_{BDE}}\overrightarrow{BD}_{\bot\overrightarrow{DE}}\quad \hbox{and}\quad \frac{l_{DE}}{2S_{CDE}}\overrightarrow{CD}_{\bot\overrightarrow{DE}}\;.
\label{e6*}
\end{equation}
On the other hand, the following equality involving a four-volume holds:
$$
24\mathcal V_{ABCDE}=\pm \det \begin{pmatrix} \overrightarrow{AD}_{\bot\overrightarrow{DE}}\\[.8ex] \overrightarrow{BD}_{\bot\overrightarrow{DE}}\\[.8ex] \overrightarrow{CD}_{\bot\overrightarrow{DE}} \end{pmatrix}\cdot l_{DE}.
$$
Thus, our desired determinant, being also the determinant composed of vectors~(\ref{e6*}), is
\begin{equation}
\det \left(\frac{\partial\vec\omega_{DE}}{\partial\psi}\right) = \pm \frac{24\mathcal V_{ABCDE}\cdot l_{DE}^2}{2S_{ADE}\cdot 2S_{BDE}\cdot 2S_{CDE}}.
\label{e7*}
\end{equation}

Multiplying (\ref{e5*}) by (\ref{e7*}) and dropping the possible minus sign which can always be eliminated by an odd basis vector permutation in some vector space, we get that {\em $\minor f_3$ has been multiplied, as a result of Pachner move $2\to 3$, by the value inverse to\/}~(\ref{tau*}). Recalling formula~(\ref{tau}) and the fact that all the change of torsion~$\tau$ is concentrated in $\minor f_3$, we see that Theorem~\ref{th_2-3} is proved.
\end{proof}

If we also recall the changes in one- and two-dimensional cells resulting from the move $2\to 3$, namely: new edge~$DE$ appeared; two-face $ABC$ disappeared, and two-faces $ADE,\allowbreak BDE,\allowbreak CDE$ appeared, we can write out the following {\em invariant\/} of moves~$2\to 3$:
\begin{equation}
\tau\cdot \frac{\prod_{\textrm{over all edges}}l^3}{\prod_{\textrm{over all 2-faces}}2S}.
\label{inv}
\end{equation}

It is already miraculous how invariant of moves~$2\leftrightarrow 3$ appears in this form of torsion multiplied by products of all volumes of certain dimensions in a simplicial complex, raised in some powers. Moreover, we will see in Subsection~\ref{ss_1-4} that formula~(\ref{inv}) gives in fact an invariant of moves~$1\leftrightarrow 4$ as well!

\subsection{Move \boldmath $1\to 4$}
\label{ss_1-4}

\begin{figure}
\begin{center}
\includegraphics[scale=1.5]{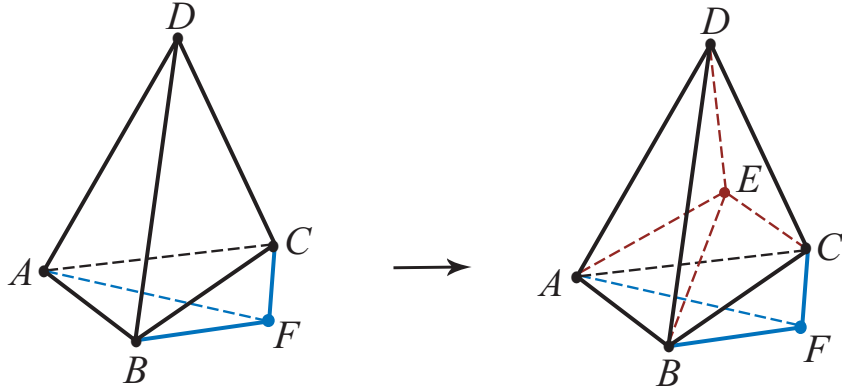}
\end{center}
\caption{Pachner move $1\to 4$ (and one adjacent tetrahedron $FABC$ which does not take part in the move)}
\label{fig2}
\end{figure}

A Pachner move $1\to 4$ is depicted in Figure~\ref{fig2}.
It consists in dividing a tetrahedron $ABCD$ into four tetrahedra $ABCE$, $BDCE$, $CDAE$ and~$ADBE$. Tetrahedron $FABC$ in Figure~\ref{fig2} does not take part in the move, but will be needed for our reasoning in this Subsection.

It will be convenient for us to represent the move $1\to 4$ as a composition of a move $0\to 2$ and the already studied move $2\to 3$. By move $0\to 2$ we understand the following change of pseudotriangulation: we slightly inflate the two-face~$ABC$ and put a new vertex~$E$ inside the resulting volume, which we then represent as consisting of two tetrahedra $ABCE$ and $EABC$ (the order of vertices reflects their mutually inverse orientation). So, these two new tetrahedra are glued to each other by their faces $AEB$, $BEC$ and~$CEA$. The face $ABC$ of tetrahedron~$ABCE$ is glued to the face of the same name of~$FABC$, while the face $ABC$ of~$EABC$ is glued to the face of the same name of~$ABCD$. As we will need to distinguish between the two resulting faces $ABC$ in the simplicial complex, we will call below the face belonging to tetrahedron $FABC$ as simply~$ABC$, while the other one --- as~$(ABC)'$. If we now do the same move $2\to 3$ in tetrahedra $EABC$ and $ABCD$ as shown in Figure~\ref{fig1}, we will get exactly the right-hand-side diagram of Figure~\ref{fig2}.

\begin{theorem}
\label{th_1-4}
The same quantity (\ref{inv}) invariant with respect to moves $2\leftrightarrow 3$ is also invariant under moves $0\to 2$ and thus under moves~$1\leftrightarrow 4$.
\end{theorem}

\begin{proof}
The following 0-, 1- and 2-cells are added to the simplicial complex as a result of move $0\to 2$:
\begin{itemize}
\item one vertex $E$,
\item three edges $AE,BE,CE$,
\item four two-faces $(ABC)',ABE,BCE,CAE$.
\end{itemize}
This leads to changes in matrices $f_2$, $f_3$ and $f_4$ in algebraic complex~(\ref{mc}). We are going to explain that, for each of these matrices, one can just replace its minor, entering in formula~(\ref{tau}), with a bigger minor, containing all `old' rows and columns as well as some `new' ones. The values of new minors are obtained from those of the old ones by simple multiplication by `small' minors, whose rows and columns correspond to the new elements in triangulation, due to the same block-triangularity reasons as in Subsection~\ref{ss_2-3}, with only a slight modification for matrix~$f_2$.
\begin{lemma}
\label{l**}
Minors in formula (\ref{tau}), taken for the initial (before move $0\to 2$) simplicial complex, can be chosen in such way that $\minor f_2$ contains the row~$d\vartheta_{ABC}$.
\end{lemma}
\begin{proof}[Proof of Lemma~\ref{l**}]
Consider a maximal nonvanishing minor of matrix~$f_2$, and let the letter~$x$ mean, for the moment, all the vertex coordinates involved in this minor (i.e., it consists of partial derivatives of $l$'s and $\vartheta$'s with respect to exactly these coordinates). Suppose the row $d\vartheta_{ABC}$ is {\em not\/} in this minor. Nevertheless, it is easy to see that the partial derivatives $\frac{\partial\vartheta_{ABC}}{\partial x}$ (substitute all the mentioned coordinates in place of~$x$!) cannot all vanish: the angle~$\vartheta_{ABC}$ does depend nontrivially on vertex coordinates. The fact that not all coordinates are in~$x$ does not prevent this: the differentials of all remaining coordinates can anyhow be made zero by means of an overall infinitesimal rotation of space~$\mathbb R^4$, which acts upon these differentials according to the mapping~$f_1$.

As $\frac{\partial\vartheta_{ABC}}{\partial x}$ is not all zero, it is possible to replace some row in $\minor f_2$ with $\frac{\partial\vartheta_{ABC}}{\partial x}$ in such way that $\minor f_2$ still does not vanish. Note that there are no problems with the corresponding change in the columns of $\minor f_3$: as we assume that complex~(\ref{mc}) is acyclic, $\minor f_3$ necessarily changes in such way that the ratio $\frac{\minor f_2}{\minor f_3}$ remains the same.
\end{proof}
We continue the proof of Theorem~\ref{th_1-4}. Assuming that the initial $\minor f_2$ has been taken as in Lemma~\ref{l**}, we add to it columns corresponding to $dw_E,dx_E,dy_E,dz_E$ and rows corresponding to $dl_{AE},dl_{BE},dl_{CE}$ and~$d\vartheta_{(ABC)'}$. Some small complication is that the row~$d\vartheta_{ABC}$ in the new minor has also to be changed with respect to the old minor. We now show how to get back to the `old'~$\vartheta_{ABC}$ while at the same time transforming $\minor f_2$ to a block triangular form. Note that $\vartheta_{ABC}$ is the dihedral angle at face~$ABC$ in four-simplex $FABCE$, while $\vartheta_{(ABC)'}$ is the angle at~$ABC$ in four-simplex $EABCD$. Thus,
$$
\vartheta_{ABC}+\vartheta_{(ABC)'}=\vartheta_{ABC}^{\rm old},
$$
where $\vartheta_{ABC}^{\rm old}$ is the angle at~$ABC$ in the four-simplex $FABCD$ {\em before\/} the move $1\to 4$. Hence, adding the row $d\vartheta_{(ABC)'}$ to the row $d\vartheta_{ABC}$ in the new minor, we get a row not depending on the coordinates of vertex~$E$, i.e., having zeros at the positions corresponding to $dw_E,dx_E,dy_E,dz_E$. Moreover, other matrix elements in this row are clearly the same as in the old minor.

The resulting minor is block-triangular, because only $dl_{AE},dl_{BE},dl_{CE}$ and~$d\vartheta_{(ABC)'}$ depend on $dw_E,dx_E,dy_E,dz_E$ (note that $\vartheta$'s for faces $AEB$, $BEC$ and $CEA$, as functions of vertex coordinates, are identical zeros). Thus,
$$
(\minor f_2)^{\rm new}=(\minor f_2)^{\rm old}\cdot \frac{dl_{AE}\wedge dl_{BE}\wedge dl_{CE}\wedge d\vartheta_{(ABC)'}}{dw_E\wedge dx_E\wedge dy_E\wedge dz_E}.
$$
To calculate the exterior form ratio, we rotate the coordinate system in such way that axis~$w$ becomes orthogonal to (the hyperplane containing) tetrahedron~$ABCE$. Then, this ratio factors in a product of values
$$
\frac{dl_{AE}\wedge dl_{BE}\wedge dl_{CE}}{dx_E\wedge dy_E\wedge dz_E}=\frac{6V_{ABCE}}{l_{AE}l_{BE}l_{CE}}
$$
(compare formulas (31) and (32) in \cite{Korepanov-JNMP}) and
$$
\frac{d\vartheta_{(ABC)'}}{dw_E}=\frac{2S_{ABC}}{6V_{ABCE}}
$$
(which is the inverse value for the height of tetrahedron $ABCE$ dropped on base $ABC$). This means that, as a result of move $0\to 2$,
\begin{equation}
\minor f_2 \quad\hbox{is multiplied by}\quad \frac{2S_{ABC}}{l_{AE}l_{BE}l_{CE}}.
\label{f2}
\end{equation}

The minor of matrix $f_3$ should be enlarged by three columns $d\vartheta_{ABE},d\vartheta_{BCE},d\vartheta_{CAE}$ and three rows corresponding to some components of the new~$\omega$, i.e., $d\vec\omega_{AE},d\vec\omega_{BE},d\vec\omega_{CE}$. Even before specifying the choice of these components, we can state the block triangularity of the new minor, caused by the fact that only $d\vartheta_{ABE},d\vartheta_{BCE},d\vartheta_{CAE}$ and no other length or angle differentials $dl$ or~$d\vartheta$ influence $d\vec\omega_{AE},d\vec\omega_{BE},d\vec\omega_{CE}$. To see this, we note that, in general, the influence on given~$d\vec\omega_i$ for a given edge~$i$ can be caused by length differentials for edges lying in the star of edge~$i$ (because the lengths of these edges determine the dihedral angles at~$i$) and by~$d\vartheta_a$ for faces~$a$ containing~$i$. As, in our situation, each of the stars of $AE,BE,CE$ consists of exactly two tetrahedra differing only in their orientations, any infinitesimal edge length changes obviously give zero $d\vec\omega$ on these edges.

To define components of $d\vec\omega$ in a convenient way, we draw three separate coordinate systems for $d\vec\omega_{AE}$, $d\vec\omega_{BE}$ and $d\vec\omega_{CE}$ respectively. One thing, however, will be common for them all, namely axis~$w$, which we draw orthogonally to (the hyperplane containing) tetrahedron~$ABCE$. The choice of one of two possible directions of~$w$, as well as other axes defined below, is not important for us. The rest of the axes must be of course orthogonal to~$w$ and to each other; we denote axis~$x$ for $d\vec\omega_{AE}$ as~$x_{AE}$ and so on, and we further specify the directions of these axes as follows:
\begin{itemize}
\item $x_{AE}$ is orthogonal to the two-dimensional plane $ABE$,
\item $y_{AE}$ lies within the plane $ABE$ and is orthogonal to $AE$,
\item $z_{AE}$ goes in or against the direction of $AE$;
\end{itemize}
\smallskip
\begin{itemize}
\item $x_{BE}$ is orthogonal to the plane $BCE$,
\item $y_{BE}$ lies within the plane $BCE$ and is orthogonal to $BE$,
\item $z_{BE}$ goes in or against the direction of $BE$;
\end{itemize}
\smallskip
\begin{itemize}
\item $x_{CE}$ is orthogonal to the plane $CAE$,
\item $y_{CE}$ lies within the plane $CAE$ and is orthogonal to $CE$,
\item $z_{CE}$ goes in or against the direction of $CE$.
\end{itemize}

So, our matrix $f_3^{\rm new}$ will contain rows corresponding to these components of~$d\vec\omega$. Consider its following submatrix:
\begin{equation}
\begin{pmatrix}
\dfrac{\partial(\omega_{AE})_{wx}}{\partial\vartheta_{ABE}} & \dfrac{\partial(\omega_{AE})_{wx}}{\partial\vartheta_{BCE}} & \dfrac{\partial(\omega_{AE})_{wx}}{\partial\vartheta_{CAE}} \\[2ex]
\dfrac{\partial(\omega_{AE})_{yw}}{\partial\vartheta_{ABE}} & \dfrac{\partial(\omega_{AE})_{yw}}{\partial\vartheta_{BCE}} & \dfrac{\partial(\omega_{AE})_{yw}}{\partial\vartheta_{CAE}} \\[2ex]
\dfrac{\partial(\omega_{BE})_{wx}}{\partial\vartheta_{ABE}} & \dfrac{\partial(\omega_{BE})_{wx}}{\partial\vartheta_{BCE}} & \dfrac{\partial(\omega_{BE})_{wx}}{\partial\vartheta_{CAE}}
\end{pmatrix}.
\label{e8*}
\end{equation}
Here, when we take the component $(\omega_{AE})_{wx}$, the letter~$x$ stands, of course, for~$x_{AE}$; similarly, $y=y_{AE}$ in $(\omega_{AE})_{yw}$, and $x=x_{BE}$ in $(\omega_{BE})_{wx}$.
\begin{lemma}
\label{l***}
Matrix (\ref{e8*}) has the following form:
\begin{equation}
\begin{pmatrix} \pm 1 & 0 & * \\ 0 & 0 & \pm\sin\varphi_{AE} \\ * & \pm 1 & 0 \end{pmatrix}
\label{e9*}
\end{equation}
and thus its determinant is $\pm \sin\varphi_{AE}$. Here $\varphi_{AE}$ is the dihedral angle at edge~$AE$ in tetrahedron~$ABCE$.
\end{lemma}
\begin{proof} Consider, for instance, the equality $\frac{\partial(\omega_{AE})_{wx}}{\partial\vartheta_{ABE}}=\pm1$. It holds because a rotation {\em around\/} face~$ABE$ goes exactly {\em in the plane\/}~$wx_{AE}$, so, the component $d(\omega_{AE})_{wx}$ is nothing but $\pm d\vartheta_{ABE}$; the sign depends on specific choices of axis directions and, as we have agreed, is not important for us, because we are calculating our invariant to within a sign. The equality $\frac{\partial(\omega_{AE})_{yw}}{\partial\vartheta_{CAE}}=\pm\sin\varphi_{AE}$ is established almost as easily. Of other equalitites for respective components of matrices (\ref{e8*}) and~(\ref{e9*}), we mention $\frac{\partial(\omega_{AE})_{wx}}{\partial\vartheta_{BCE}}=0$ which holds for a different reason: $\omega_{AE}$ simply does not depend on the angle at face~$BCE$, because this latter does not contain the edge~$AE$.
\end{proof}
Recall that we are still within the proof of Theorem~\ref{th_1-4}. Joining the result of Lemma~\ref{l***} with block triangularity argument, we can state that
\begin{equation}
\minor f_3 \quad \hbox{is multiplied by} \quad \sin\varphi_{AE}=\frac{6V_{ABCE}\cdot l_{AE}}{2S_{ABE}\cdot 2S_{CAE}}.
\label{f3}
\end{equation}
It remains to study the change in $\minor f_4$. Here, six columns are added, corresponding to the six not-yet-used components of $d\vec\omega_{AE}$, $d\vec\omega_{BE}$ and~$d\vec\omega_{CE}$, and six rows corresponding to the six-component quantity~$d\rho_E$. Block triangularity is ensured by the fact that $d\rho_E$ is influenced only by $d\vec\omega$ of those edges which end in~$E$. So, $(\minor f_4)^{\rm old}$ is multiplied by a $6\times 6$ determinant in order to obtain $(\minor f_4)^{\rm new}$ and, besides, this determinant factorizes in two $3\times 3$ determinants, the first of which corresponds to components of $d\vec\omega$ and~$d\rho$ orthogonal to axis~$w$, while the second --- to components of type $wx$ or~$yw$.

The first $3\times 3$ determinant can be written as follows:
\begin{equation}
\frac{d(\omega_{AE})_{xy}\wedge d(\omega_{BE})_{xy}\wedge d(\omega_{CE})_{xy}}{\bigwedge d\rho_{\perp w}}.
\label{e10*}
\end{equation}
Remember that here `$xy$' are different for different~$d\omega$. It is not hard to understand that (\ref{e10*}) equals the volume of parallelepiped built on unit vectors going in the directions of axes $z_{AE},z_{BE},z_{CE}$, that is,
\begin{equation}
\frac{6V_{ABCE}}{l_{AE}l_{BE}l_{CE}}.
\label{f4*}
\end{equation}
The second $3\times 3$ determinant can be written as
$$
\frac{d(\omega_{BE})_{yw}\wedge d(\omega_{CE})_{wx}\wedge d(\omega_{CE})_{yw}}{\bigwedge d\rho_{\| w}}.
$$
It equals the volume of parallelepiped built on unit vectors going along axes $y_{BE}$, $x_{CE}$ and~$y_{CE}$. It is not hard to see that this volume equals the sine of angle between edges $BE$ and~$CE$, which is
\begin{equation}
\frac{2S_{BCE}}{l_{BE}l_{CE}}.
\label{f4**}
\end{equation}
Combining the factors (\ref{f4*}) and~(\ref{f4**}), we can state that
\begin{equation}
\minor f_4 \quad \hbox{is multiplied by} \quad \frac{6V_{ABCE}\cdot 2S_{BCE}}{l_{AE}l_{BE}^2 l_{CE}^2}.
\label{f4}
\end{equation}

And now, combining (\ref{f2}), (\ref{f3}) and~(\ref{f4}), we see that torsion $\tau$, defined according to formula~(\ref{tau}), is multiplied under the move $0\to 2$ by
$$
\frac{\prod_{\textrm{over new faces}} 2S}{\prod_{\textrm{over new edges}} l^3}.
$$
This means that Theorem \ref{th_1-4} has been proved.
\end{proof}

\section{Discussion of results}

In this paper, we have not given concrete examples of calculations. In fact, we plan to use our complex~(\ref{mc}) not so much in its form presented in this paper as the basis for further modifications. In the forthcoming paper~\cite{in_progress} we intend, first, to demonstrate interesting properties of reducibility (possibility to be decomposed into a direct sum) for complex~(\ref{mc}) and, second, make concrete calculations for some framed knots, using a suitable modification of complex~(\ref{mc}) made in analogy with paper~\cite{DKM}.

Of other ideas for the nearest future, we mention the idea of geometrizing simplicial complexes by putting more complicated geometric objects in correspondence to their elements. For instance, not only just a point in a Euclidean or other space can correspond to a vertex in a simplicial complex but, say, a line segment lying in that space. We understand though that any object can be thought of as a point in a space of such objects. It is worth thinking how to choose such a space so that it would deserve to be called quantum.

\subsection*{Acknowledgments} Problems of quantization were one of the favorite themes of my scientific supervisor F.A.~Berezin. He also drew my attention to ideas related to integrable models in mathematical physics. I would like to use this opportunity to remember Felix Alexandrovich with gratitude.

This paper has been written with a partial financial support from Russian Foundation for Basic Research, Grant no.~04-01-96010.

\end{document}